\documentclass[12pt]{amsart}

\usepackage{color,graphicx}
\usepackage{psfrag}
\usepackage{amssymb}                    
\usepackage{amsmath}                    
\usepackage{amsthm}                     
\usepackage{latexsym}                   
\usepackage{subfigure}
\usepackage{epsfig}                     
\usepackage{pstricks,pst-node,pst-text,pst-tree}
\usepackage{amscd}
\usepackage{amsmath}
\usepackage{psfrag}

\newtheorem{theorem}{Theorem}[section]
\newtheorem{corollary}[theorem]{Corollary}
\newtheorem{lemma}[theorem]{Lemma}

\newtheorem{proposition}[theorem]{Proposition}

\theoremstyle{definition}
\newtheorem{definition}[theorem]{Definition}

 \newtheorem{remark}[theorem]{Remark}

\DeclareMathOperator{\Spec}{Spec}
\DeclareMathOperator{\Proj}{Proj}

\DeclareMathOperator{\face}{face}

\newcommand{\one}{\ensuremath{(\mathrm{i})}}
\newcommand{\two}{\ensuremath{(\mathrm{ii})}}

 \newcommand{\Chow}{\ensuremath{\operatorname{Chow}(X_P/T_L)}}

 \newcommand{\K}{\Bbbk}

 \newcommand{\N}{\ensuremath{\mathbb{N}}}

 \newcommand{\Q}{\ensuremath{\mathbb{Q}}}

 \newcommand{\Z}{\ensuremath{\mathbb{Z}}}

 \newcommand{\git}{\ensuremath{\operatorname{/\!\!/}}}

 \newcommand{\image}{\ensuremath{\operatorname{im}_{\mathbb Z}}}
 
 \newcommand{\owe}{\ensuremath{\mathcal{O}}}
 
 \newcommand{\poset}{\ensuremath{\mathcal{P}}}
 \newcommand{\rec}{\ensuremath{\operatorname{rec}}}
 \newcommand{\relint}{\ensuremath{\operatorname{relint}}}
 \newcommand{\sat}{\ensuremath{\operatorname{sat}}}

 \numberwithin{equation}{section}
\begin{document}

 \bibliographystyle{amsalpha} 
 
 \title{Fiber fans and toric quotients}
 
 \author{Alastair Craw} \address{Department of Mathematics, Stony
   Brook University, Stony Brook, NY 11794}
 \email{craw@math.sunysb.edu}
 
 \author{Diane Maclagan} \address{Department of Mathematics,
       Hill Center-Busch Campus,
       Rutgers University,
       110 Frelinghuysen Rd,
       Piscataway, NJ 08854} \email{maclagan@math.rutgers.edu}

\date{\today}

 \begin{abstract}
   The GIT chamber decomposition arising from a subtorus action on a
   quasiprojective toric variety is a polyhedral complex.  Denote by
   $\Sigma$ the fan that is the cone over the polyhedral complex.
   In this paper we show that the toric variety defined by the fan
   $\Sigma$ is the normalization of the toric Chow quotient of a
   closely related affine toric variety by a complementary torus.
 \end{abstract}

 \maketitle

\section{Introduction}

Let $X_P$ be the quasiprojective toric variety defined by a
full-dimensional polyhedron $P\subseteq \Q^n$, and let $T_L = \Spec
\K[L]$ be a $d$-dimensional subtorus of the dense torus $(\K^*)^n$ of
$X_P$.  The inclusion of $T_L$ into $(\K^*)^n$ induces a linear map
$\pi\colon \Q^n\rightarrow \Q^d \cong L\otimes_{\Z} \Q$ and hence a
linear projection of polyhedra $\pi\colon P\rightarrow Q$, where $Q =
\pi(P)$.  Thaddeus~\cite{Thaddeus} proved that $\pi$ induces a cell
decomposition of $Q$ that is the GIT chamber decomposition arising
from the action of $T_L$ on $X_P$.  See Hu~\cite{Hu} for the most
general results about toric GIT, including quotients of toric
varieties that do not necessarily arise from polyhedra.

By taking the cone over each cell in this GIT chamber decomposition,
we obtain a fan $\Sigma$ in $\Q^d\oplus \Q$.  It is natural to ask for
an explicit description of the toric variety defined by this fan.  The
main result of this paper provides an answer to this question as follows (see
Theorem~\ref{thm:main} and Section~\ref{sec:GeneralDuality} for relevant notation):

\begin{theorem}
  \label{t:main}
  Let $T_L$ be a subtorus acting on a toric variety $X_P$, and let
  $\Sigma$ be the fan arising from the GIT chamber decomposition
  associated to the $T_L$-action on $X_P$ as above.  Then there is a
  polyhedral cone $P':= (\widetilde{P})^\vee$ and a lattice
  $L':=\image(\widetilde{\pi}^\vee)$ such that $\Sigma$ is the fan of
  the toric Chow quotient $X_{P'}/T_{L'}$.
\end{theorem}

The Chow quotient $X_{P'}/T_{L'}$ of a \emph{projective} toric variety
$X_{P'}$ by a subtorus $T_{L'}$ was studied by
Kapranov--Sturmfels--Zelevinsky~\cite{KSZ}.  We generalize this by
constructing the \emph{toric Chow quotient} $X_{P}/T_{L}$ of an
arbitrary quasiprojective toric variety $X_{P}$ by a subtorus $T_{L}$
(see Definition~\ref{defn:toricchow}).

In order to prove Theorem~\ref{t:main}, we consider the \emph{fiber
  fan} associated to a linear projection of polyhedra $\pi\colon
P\rightarrow Q$, denoted $\mathcal{N}(P,Q)$, generalizing the normal
fan of a fiber polytope (see
Billera--Sturmfels~\cite{BilleraSturmfels}).  We show that the toric
Chow quotient $X_P/T_L$ is a not-necessarily normal toric variety
whose fan is the fiber fan $\mathcal{N}(P,Q)$.  In the special case
where $P$ is a polytope, so $X_P$ is projective, this result is due to
Kapranov--Sturmfels--Zelevinsky~\cite{KSZ}.  Note however that the
statement of Theorem~\ref{t:main} is new even in the case where $P$ is
a polytope.

If the original polyhedron $P$ is a polyhedral
cone, the statement of Theorem~\ref{t:main} can be streamlined, and a
converse can be added as follows (see Theorem~\ref{thm:affineduality}).
  
 \begin{theorem}[Duality for affine toric quotients]
 \label{t:maintwo}
   Let $P\subseteq \Q^n$ be a polyhedral cone and $\pi\colon
   P\rightarrow Q$ a linear projection.  Set $P' := P^\vee$ and $L':=
   \image(\pi^\vee)$.  Then
 \begin{enumerate}
 \item[\one] the fan of the toric Chow quotient $X_{P}/T_{L}$ is equal
   to the GIT chamber decomposition arising from
   the action of $T_{L'}$ on $X_{P'}$; and
 \item[\two] the GIT chamber decomposition arising from the
   $T_L$-action on $X_P$ is the fan of the toric Chow quotient
   $X_{P'}/T_{L'}$.
 \end{enumerate}
 \end{theorem}
 
 The statement of Theorem~\ref{t:maintwo}, part \one, was partially
 known to Altmann--Hausen~\cite{AltmannHausen} in the case where
 $P$ is a simplicial cone, though no proof was given.
 
 The main theorems are proved in Section~\ref{s:mainproof}.
 Section~\ref{s:polyhedral} contains the necessary polyhedral
 preliminaries, while Sections~\ref{sec:toricGIT} and
 \ref{sec:toricChow} contain background on toric GIT, and introduce
 the toric Chow quotient, respectively.

 \medskip

 \noindent \textbf{Acknowledgements} Special thanks are due to Rekha
 Thomas for providing a great deal of help.  This project arose out of
 our joint work in \cite{CMT1}.  We also thank Klaus Altmann and Mark
 Haiman for useful conversations.  The second author was partially
 supported by NSF grant DMS-0500386.

\section{Polyhedral Geometry}
 \label{s:polyhedral}
\subsection{Polyhedral Conventions}

Let $P \subseteq \mathbb Q^n$ be an $n$-dimensional polyhedron.  For
$\mathbf{w} \in (\mathbb Q^n)^*$ we denote by $\face_{\mathbf{w}}(P)$
the face of $P$ minimizing $\mathbf{w}$.  Given a face $F$ of $P$, the
\emph{inner normal cone} $\mathcal N_{P}(F)$ is the set of ${\textbf
  w} \in (\Q^n)^{\ast}$ such that $\face_{\mathbf{w}}(P)=F$.  The
\emph{inner normal fan} $\mathcal{N}(P)$ of $P$ is the polyhedral fan
whose cells are the inner normal cones $\{\mathcal N_{P}(F)\}$ as $F$
varies over the faces of $P$.  Two polyhedra $P$ and $P'$ are {\em
  normally equivalent} if they have the same normal fan.

The {\em recession cone} $\rec(P)$ of a polyhedron $P$ is $\{ {\bf
  u}\in \Q^n : {\bf x} + {\bf u} \in P \text{ for all } {\bf x} \in P
\}$.  If $P=\{ {\bf x} \in \mathbb Q^n : A{\bf x} \geq {\bf b} \}$ for
some $r \times n$ matrix $A$ and vector ${\bf b} \in \mathbb Q^r$,
then $\rec(P)=\{ {\bf x} \in \mathbb Q^n : A{\bf x} \geq 0 \}$.  The
{\em support} of a polyhedral fan is the set of vectors lying in some
cone of the fan.  The fan $\mathcal{N}(P)$ is supported on the dual
cone $\rec(P)^\vee = \{ {\bf y} \in (\Q^k)^* : {\bf y} \cdot {\bf x}
\geq 0 \,\, \forall \,\, {\bf x} \in \rec(P) \}$ of the
recession cone $\rec(P)$ of $P$.  Indeed, ${\bf y}\not\in
\rec(P)^\vee$ if and only if there exists ${\bf x}\in \rec(P)$ such
that ${\bf y}\cdot {\bf x}< 0$, in which case ${\bf y}\cdot {\bf x}'$
is unbounded below for ${\bf x}':= \lambda {\bf x}\in \rec(P)$ with
$\lambda \in \Q_{\geq 0}$.

Let $\pi \colon \mathbb Q^n \rightarrow \mathbb Q^d$ be a linear map,
and set $Q:=\pi(P)$.  We may assume that $\pi$ is the projection onto
the last $d$ coordinates.  We denote by $\pi_{n-d}$ the projection
onto the first $n-d$ coordinates $\Q^{n-d}$, 
and for $\mathbf{q} \in Q$ we set $P_\mathbf{q}:=
\pi_{n-d}(\pi^{-1}(\mathbf{q}) \cap P)$.
For ${\bf w}\in (\Q^{n-d})^*$, we write
$\widetilde{\face}_{\bf{w}}(P_\mathbf{q})$ for the unique smallest
face of $P$ containing the preimage of
$\face_{\mathbf{w}}(P_{\mathbf{q}})$ under $\pi_{n-d}$.

\subsection{Fiber fans}
In this section we define the notion of a fiber fan $\mathcal N(P,Q)$
for a linear projection $\pi: P \rightarrow Q$ of polyhedra.  This is
the common refinement of the normal fans of the fibers
$\{\pi^{-1}({\bf q}): {\bf q} \in Q\}$ of the map $\pi$.  We first
show that this definition makes sense.

\begin{lemma}
For $\mathbf{q}, \mathbf{q}' \in Q$, set $\mathbf{q} \sim
  \mathbf{q}$ if $P_\mathbf{q}$ is normally equivalent to
  $P_{\mathbf{q}'}$.  Then there are only finitely many equivalence
  classes for $\sim$.
\end{lemma}

\begin{proof}
  Write $P= \{ {\bf x} : A{\bf x} \geq {\bf b} \}$, where $A$ is an $r
  \times n$ matrix for some $r$, and ${\bf b} \in \mathbb Q^r$.  Then
  for $\mathbf{q} \in Q$, $\pi^{-1}(\mathbf{q}) \cap P = \{ ({\bf
    x}',\mathbf{q}) : {\bf x}' \in \mathbb Q^{n-d}, A({\bf x}',{\bf
    q}) \geq {\bf b} \}$.  Write $A$ in block form as $(A'|A'')$,
  where $A'$ is an $r \times (n-d)$ matrix, and $A''$ is an $r \times
  d$ matrix.  Then $P_{\mathbf{q}}= \{ {\bf x}' \in
  \mathbb Q^{n-d} : A' {\bf x}' \geq {\bf b}-A''\mathbf{q}\}$.  Now
  for a given matrix $A'$ there are only finitely many normal
  equivalence classes of polyhedra of the form $\{{\bf x}'\in \mathbb
  Q^{n-d} : A'{\bf x}' \geq {\bf c} \}$ as $\mathbf{c}$ varies, so it
  follows that there are only finitely many normal equivalence classes
  of fibers of $\pi$.
\end{proof}

\begin{definition}
  Pick $\mathbf{q}_1,\dots,\mathbf{q}_r \in Q$ such that the
  $P_{\mathbf{q}_i}$ are representatives of the different normal
  equivalence classes of fibers of $\pi$.  Let $F$ be the Minkowski
  sum of the $P_{\mathbf{q}_i}$.  The {\em fiber fan} $\mathcal
  N(P,Q)$ is the inner normal fan of $F$.
\end{definition}

\begin{lemma}
The recession cone of $P_\mathbf{q}$ is the same for all $\mathbf{q}
  \in Q$.  
\end{lemma}

\begin{proof}
  Write $P=\{{\bf x} \in \mathbb Q^n : A{\bf x} \geq {\bf b} \}$ for
  some $r \times n$ matrix $A$ and ${\bf b} \in \mathbb Q^r$.
  Decompose $A=(A'|A'')$, where $A'$ is an $r \times (n-d)$ matrix,
  and $A''$ is an $r \times d$ matrix.  Translating the polyhedron $P$
  by ${\bf a} \in \mathbb Q^n$ translates $Q$ by $\pi({\bf a})$, and
  $P_\mathbf{q} = (P+{\bf a})_{\mathbf{q}+\pi({\bf a})}$, so we may
  assume ${\bf 0}\in P$.  
Note that $P_{\bf 0} = \{ {\bf x} \in \mathbb Q^{n-d} : A' {\bf x} \geq
  {\bf b} \}$. We will show that $\rec(P_{\mathbf{q}}) = \rec(P_{\bf
    0})$ for all ${\bf q}\in Q$.
  
  Let $\mathbf{v} \in \rec(P_{\mathbf{q}})$.  Then ${\bf u}+\lambda
  \mathbf{v} \in P_\mathbf{q}$ for all ${\bf u} \in P_\mathbf{q}$ and
  $\lambda>0$.  Since ${\bf 0} \in P$, we have $b_i\leq 0$ for all
  $i$, where $b_i$ is the $i$th component of ${\bf b}$.  Suppose that
  $(\mathbf{v},{\bf 0}) \not \in P$ for ${\bf 0}\in \Q^d$.  Then there
  is some $i$ with ${\bf a_i} \cdot (\mathbf{v},{\bf 0}) < b_i\leq 0$,
  where ${\bf a_i}$ is the $i$th row of $A$.  But then ${\bf a_i}
  \cdot ({\bf u}+\lambda \mathbf{v}, \mathbf{q})={\bf a_i} \cdot ({\bf
    u},\mathbf{q}) + \lambda {\bf a_i} \cdot(\mathbf{v},0) < b_i$ for
  $\lambda$ sufficiently large, which means that $({\bf u}+\lambda
  \mathbf{v},\mathbf{q}) \not \in P$, so ${\bf u} +\lambda \mathbf{v}
  \not \in P_\mathbf{q}$.  Therefore we have $\{ (\mathbf{v},{\bf 0})
  \in \Q^n: \mathbf{v} \in \rec(P_{\mathbf{q}}) \} \subseteq P$ after
  all.  The above argument shows that ${\bf a_i} \cdot
  (\mathbf{v},{\bf 0}) = {\bf a'_i} \cdot \mathbf{v} \geq 0$ for each
  row ${\bf a'_i}$ of $A'$, so $\mathbf{v}\in\rec(P_0)$.
  
  For the opposite inclusion, note that since $P_{\bf 0} = \{ {\bf x}
  \in \mathbb Q^{n-d} : A' {\bf x} \geq {\bf b} \}$, the set $\{
  (\mathbf{v},0) : \mathbf{v} \in \rec(P_0) \}$ lies in $\rec(P)$.
  Thus if $\mathbf{v} \in \rec(P_0)$, then $({\bf
    u},\mathbf{q})+\lambda(\mathbf{v},{\bf 0}) \in P$ for all $\lambda
  >0$ and ${\bf u} \in P_\mathbf{q}$.  This means $({\bf u}+\lambda
  \mathbf{v},\mathbf{q}) \in P$ for all $\lambda>0$ and ${\bf u} \in
  P_\mathbf{q}$, so ${\bf u}+\lambda \mathbf{v} \in P_\mathbf{q}$ for
  all $\lambda>0$ and ${\bf u} \in P_\mathbf{q}$, giving $\mathbf{v} \in
  \rec(P_{\mathbf{q}})$ as required.
\end{proof}

If $\mathcal F_1,\dots,\mathcal F_r$ are polyhedral fans with the same
support, then their {\em common refinement} is the fan $\mathcal F$
whose cones are the intersections $\cap_{i=1}^r \sigma_i$, where
$\sigma_i$ is a cone in $\mathcal F_i$.  If $P$ and $P'$ are two
polyhedra with the same recession cone, then the normal fan of the
Minkowski sum of $P$ and $P'$ is the common refinement of the normal
fans of $P$ and $P'$. Thus the fiber fan $\mathcal N(P,Q)$ is the
common refinement of the normal fans of the fibers.

\begin{remark}
  In the case that $P$ is a polytope, the fiber fan $\mathcal N (P,
  Q)$ is the normal fan of the {\em fiber polytope} introduced by
  Billera and Sturmfels \cite{BilleraSturmfels}.  While the fiber fan
  of a linear projection of polyhedra is still the normal fan of a
  polyhedron, there is not a canonical choice of a polyhedron with
  that normal fan if the polyhedron being projected is not a polytope.
  This motivates our definition of the fiber fan, instead of a more
  general fiber polyhedron.  Note that for applications to toric
  varieties, a fan suffices to define the variety.
\end{remark}

\subsection{Duality for polyhedral cones}

 \label{sec:coneduality} 
 This section establishes a duality result for polyhedral cones that
 is the key to Theorem~\ref{t:maintwo}. Versions of
 Lemmas~\ref{lemma:faceequality} and \ref{l:duallemma} are known to
 the experts, though we know of no proofs in the literature.

 Let $C\subset \Q^n$ be a full-dimensional polyhedral cone.  For $d <
 n$, let $\pi\colon \Q^n\rightarrow \Q^d$ be the projection onto the
 last $d$ coordinates.  For $\mathbf{v}\in \pi(C)$, consider the
 polyhedral slice $C_{\mathbf{v}}\subseteq \mathbb Q^{n-d}$.  As
 $\mathbb Q^{n-d}$ is canonically isomorphic to $\ker(\pi)$, we may
 regard the normal fan $\mathcal{N}(C_\mathbf{v})$ as lying in the
 dual vector space $\ker(\pi)^* \cong (\mathbb Q^{n-d})^*$.  Let
 $\pi^\vee\colon (\Q^n)^*\rightarrow \ker(\pi)^*$ denote the map dual
 to the inclusion of $\ker(\pi)$ in $\Q^n$.  In the basis of
 $(\Q^n)^*$ dual to the standard basis of $\Q^n$ this is projection on
 the first $n-d$ coordinates.  For ${\textbf w}\in \pi^{\vee}(C^\vee)$
 we consider the polyhedral slice $C^{\vee} \cap
 (\pi^{\vee})^{-1}(\textbf w)$, and write $C^\vee_{\bf w}:=
 \pi_d(C^\vee\cap (\pi^{\vee})^{-1}({\bf w})) \subseteq
(\Q^d)^*$ 
for the isomorphic image of $C^{\vee} \cap (\pi^{\vee})^{-1}(\textbf
w)$ under the projection $\pi_d\colon (\Q^n)^*\rightarrow
(\Q^d)^*$ onto the last $d$ coordinates.
There is a canonical isomorphism $\ker(\pi^\vee)\cong (\Q^d)^*$, so
the normal fan $\mathcal{N}(C^\vee_{\textbf w})$ can be regarded as
living in either $\Q^d$ or in $(\ker(\pi^{\vee}))^*$ for all ${\bf
    w}\in (\Q^{n-d})^*$.

Let $\relint(F)$ denote the relative interior of a set $F$, which is
the interior of $F$ in its affine span.

 \begin{lemma} 
 \label{lemma:faceequality} 
 Let $F\subset C$ be a face.   If $\mathbf{v} \in \pi(\relint(F))$
 then we have
 $\pi^{\vee}( \mathcal{N}_C(F)) =
 \mathcal{N}_{C_\mathbf{v}}(F_\mathbf{v})$. 
 \end{lemma}

 \begin{proof}
   Write $C = \{{\bf x}\in \Q^n : A{\bf x}\geq {\bf 0}\}$ where
   $A$ is the $r\times n$ matrix whose rows  ${\bf a}_1,\dots
   ,{\bf a}_r\in (\Q^n)^*$ form the facet normals of $C$.  We may
   assume that $F$ is cut out by the first $p$ facet inequalities
   defining $C$, so $\mathcal{N}_{C}(F) = \Q_{\geq 0}\langle{\bf
     a}_1, \ldots, {\bf a}_p\rangle$ and
   $\pi^{\vee}(\mathcal{N}_{C}(F)) = \Q_{\geq 0}\langle
   \pi^{\vee}({\bf a}_1), \dots, \pi^{\vee}({\bf a}_p)\rangle$.  Note that
 \[
 F\cap \pi^{-1}(\mathbf{v}) = \Big{\{} ({\bf u},\mathbf{v})  :
 \begin{array}{ll} {\bf a}_i \cdot ({\bf u} , \mathbf{v}) = 0 &
     \text{for } i = 1, \ldots, p; \\ {\bf a}_i \cdot ({\bf u}, \mathbf{v})
     \geq  0 & \text{for } i = p+1, \ldots, r. \end{array}\Big{\}}
 \]
 since $F$ is cut out by the first $p$ facet inequalities defining
 $C$.  The assumption that $\mathbf{v} \in \pi(\relint(F))$ ensures
 that each $\geq$ above is a $>$, so 
 \[
 F_{\mathbf{v}} = \big{\{}{\bf u} \in \Q^{n-d} : \pi^{\vee}({\bf a}_i)
 \cdot {\bf u} = c_i  \text{ for }i = 1, \dots,p; \;
 \pi^{\vee}({\bf a}_i) \cdot {\bf u} > c_i  \text{ for }i =
 p+1,\dots, r \big{\}}
 \]
 where $c_i=-\mathbf{a}_i \cdot \mathbf{v}$ for $i=1, \ldots, r$.
 This implies that the rational cone
 $\mathcal{N}_{C_{\mathbf{v}}}(F_\mathbf{v})$ is also generated by
 $\{\pi^{\vee}({\bf a}_1), \dots, \pi^{\vee}({\bf a}_p)\}$ as claimed.
 \end{proof}

 \begin{lemma} 
 \label{l:duallemma} 
 Fix $\mathbf{v}\in \pi(C)$ and ${\bf w}\in \pi^{\vee}(C^\vee)$.
 If $F\subset C$ is a face then
 \[
 F = \widetilde{\face}_{\bf w}(C_{\mathbf{v}})\iff \mathbf{v}
 \in\pi(\relint(F))\text{ and }{\bf w}
 \in\pi^{\vee}\big{(}\relint(\mathcal{N}_{C}(F))\big{)}.
 \] 
 \end{lemma}

 \begin{proof}
   Suppose $F = \widetilde{\face}_{\bf w}(C_{\mathbf{v}})$. Then the
   preimage of $\face_{\bf w}(C_{\mathbf{v}})$ under $\pi_{n-d}$ lies
   in no proper face of $F$, so the intersection of $\relint(F)$ with
   this preimage is nonempty.  This implies the set
   $\pi^{-1}(\mathbf{v}) \cap \relint(F)$ is nonempty, so $\mathbf{v}
   \in \pi(\relint(F))$.  We have ${\bf w}\in
   \relint(\mathcal{N}_{C_\mathbf{v}}(\face_{\bf w}(C_\mathbf{v})))$
   by definition.  Applying Lemma~\ref{lemma:faceequality} to the face
   $F_\mathbf{v} = \pi_n(F\cap \pi^{-1}(\mathbf{v})) = \face_{\bf
     w}(C_\mathbf{v})$ gives ${\bf w}\in
   \relint\pi^{\vee}(\mathcal{N}_C(F)) =
   \pi^{\vee}(\relint\mathcal{N}_C(F))$.
   
   Conversely, suppose $\mathbf{v} \in \pi(\relint(F))$ and ${\bf w}
   \in\pi^{\vee}(\relint(\mathcal{N}_{C}(F)))$.  The first assumption
   ensures that the face $F_\mathbf{v}\subseteq C_\mathbf{v}$ is nonempty
   and its preimage under $\pi_{n-d}$ lies in no proper subface of $F$.
   Also, ${\bf w} \in\pi^{\vee}(\relint(\mathcal{N}_{C}(F))) =
   \relint\pi^{\vee}(\mathcal{N}_C(F)) = \relint
   \mathcal{N}_{C_\mathbf{v}}(F_\mathbf{v})$ by
   Lemma~\ref{lemma:faceequality}, hence $F_\mathbf{v} = \face_{\bf
     w}(C_{\mathbf{v}})$.  Since the preimage of $F_\mathbf{v}$ under
   $\pi_{n-d}$ lies in no proper subface of $F$, we have $F =
   \widetilde{\face}_{\bf w}(C_{\mathbf{v}})$.
 \end{proof}

\begin{figure}[!ht]
 \input{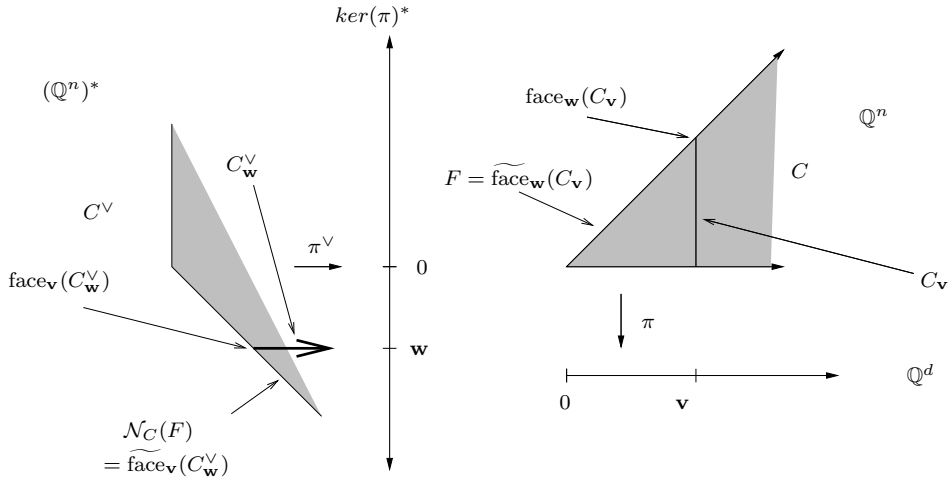}
 \caption{Figure for Proposition~\ref{p:fiberduality}}
 \label{theorem2.3}
 \end{figure}

 We now present the main result of this section.  See
 Figure~\ref{theorem2.3} for an illustration.

 \begin{proposition}
 \label{p:fiberduality} 
 For $\mathbf{v}\in \pi(C)$ and ${\bf w}\in \pi^{\vee}(C^\vee)$ we
 have
 \begin{equation*}
 \label{eqn:fiberduality}
 \mathcal{N}_C\big{(}\widetilde{\face}_{\bf w}(C_{\mathbf{v}})\big{)} =
 \widetilde{\face}_{\mathbf{v}}(C^\vee_{\bf w}).
 \end{equation*}
 \end{proposition}
 \begin{proof}
   Applying Lemma~\ref{l:duallemma} to the face
   $F:=\widetilde{\face}_{\bf w}(C_{\mathbf{v}})$ of $C$ gives
   $\mathbf{v} \in \pi(\relint (F))$ and ${\bf w}\in
   \pi^{\vee}(\relint(\mathcal{N}_{C}(F)))$. The face
   $\mathcal{N}_C(F)$ of $C^{\vee}$ satisfies
   $\mathcal{N}_{C^{\vee}}(\mathcal{N}_C(F)) = F$, so ${\bf w} \in
   \pi^{\vee}(\relint(\mathcal{N}_C(F)))$ and $\mathbf{v} \in
   \pi(\relint(\mathcal{N}_{C^{\vee}}(\mathcal{N}_C(F))))$.  Now
   apply Lemma~\ref{l:duallemma} to the face $\mathcal{N}_C(F)$ of
   $C^{\vee}$ to deduce that $\mathcal{N}_C(F) =
   \widetilde{\face}_{\mathbf{v}}(C^\vee_{\bf w})$ as required.
 \end{proof}

 \begin{remark}
   Proposition~\ref{p:fiberduality} is the key duality result that
   is required in the proof of Theorem 6.2 of
   Craw--Maclagan--Thomas~\cite{CMT1}.
 \end{remark}

 The following definition is a variant for polyhedral cones of a
 notion for polytopes due to Billera--Sturmfels~\cite[\S
 2]{BilleraSturmfels}.

 \begin{definition}
   For $\mathbf{v} \in \pi(C)$, the subdivision of
   ${\pi^\vee}(C^\vee)$ consisting of the cells $\big\{ \pi^\vee\big(
   \widetilde{\face}_{\mathbf{v}} (C^\vee_{\bf w})\big{)} : {\bf w}
   \in {\pi^\vee}(C)\big\}$ is the ${\mathbf{v}}$\emph{-induced
     ${\pi^\vee}$-coherent subdivision} of ${\pi^\vee}(C^\vee)$.
 \end{definition} 

 The duality result of Proposition~\ref{p:fiberduality} leads
 immediately to the following description of the normal fan of the
 polyhedral slice $C_\mathbf{v}$ for $\mathbf{v}\in \pi(C)$.

\begin{corollary} 
\label{cor:normalfan}
For $\mathbf{v}\in \pi(C)$, the inner normal fan of the polyhedron
$C_{\mathbf{v}}$ is the $\mathbf{v}$-induced $\pi^{\vee}$-coherent
subdivision of $\pi^{\vee}(C^\vee)$. In other words,
 \[
 \mathcal{N}(C_\mathbf{v}) =
 \big{\{}\pi^{\vee}\big{(}\widetilde{\face}_{\mathbf{v}}(C^\vee_{\bf
 w})\big{)} : {\bf w} \in \pi^{\vee}(C^\vee) \big{\}}.
 \]  In
 particular, for $\mathbf{v} = \emph{\textbf{0}}$ we have
 $\mathcal{N}(C_{\emph{\textbf{0}}}) = \pi^{\vee}(C^\vee)$.
\end{corollary}

\begin{proof}
  By Lemma~\ref{lemma:faceequality} and
  Proposition~\ref{p:fiberduality}, we have
  $$\mathcal{N}_{C_{\mathbf{v}}}\big{(}\face_{\bf w}(C_{\mathbf{v}})\big{)} =
  \pi^{\vee}\big{(}\mathcal{N}_{C}\big{(}\widetilde{\face}_{\bf
  w}(C_{\mathbf{v}})\big{)}\big{)} =
  \pi^{\vee}\big{(}\widetilde{\face}_{\mathbf{v}}(C^\vee_{\bf
  w})\big{)}.$$ The set of cones $\mathcal{N}_{C_{\mathbf{v}}}(\face_{\bf
  w}(C_{\mathbf{v}}))$ taken over all vectors ${\bf w} \in
  \pi^{\vee}(C^\vee)$ equals the inner normal fan
  $\mathcal{N}(C_\mathbf{v})$ as required.  For the second statement,
  observe that $\face_{\textbf{0}}(C^\vee_{\bf w}) = C^\vee_{\bf w}$
  for all ${\bf w} \in \pi^{\vee}(C^\vee)$, and the smallest face of
  the cone $C^\vee$ containing the set $C^\vee_{\bf w}$ is
  $C^\vee$ itself.
\end{proof}

 \begin{definition}
 \label{d:subdivisioninduced}
 Given a linear map $\pi\colon \Q^n\rightarrow \Q^d$ and a fan
 $\Sigma\subseteq \Q^n$, the \emph{fan induced by $\pi(\Sigma)$} is
 the fan obtained by taking the common refinement of the cones
 $\{\pi(\sigma) : \sigma\in \Sigma\}$.  If $\Sigma$ consists of the
 faces of a cone $P\subseteq \Q^n$, we simply write $\pi(P)$.
\end{definition}

\begin{corollary}
 \label{cor:coneduality}
The fiber fan $\mathcal N(C,\pi(C))$ is the fan induced by
$\pi^{\vee}(C^{\vee})$.
\end{corollary}
\begin{proof}
The fan $\mathcal N(C, \pi(C))$ is the common refinement of the
normal fans of the fibers $C_{\mathbf{v}}$ for $\mathbf{v} \in
\pi(C)$.  As $\mathbf{v}$ varies in $\pi(C)$,
$\widetilde{\face}_{\mathbf{w}}(C^{\vee}_{\mathbf{v}})$ varies over
all faces of $C^{\vee}$.  Thus $\mathcal N(C,\pi(C))$ is the common
refinement of all $\pi^{\vee}(F)$ for $F$ a face of $C^{\vee}$.
\end{proof}


 \subsection{Duality for general polyhedra}
 \label{sec:GeneralDuality}
 We now generalize the results of the previous subsection from
 polyhedral cones to general polyhedra.

 \begin{definition}
 \label{d:tildeP}
 For any polyhedron $P \subseteq \mathbb Q^n$, let $\widetilde{P}
 \subseteq \Q^n\oplus \Q$ be the polyhedral cone obtained as the
 closure in $\Q^{n}\oplus \Q$ of the set $\{ (\lambda {\bf x},\lambda)
 : {\bf x} \in P, \lambda \in \Q_{>0} \}$.  Note the $\widetilde{P} \cap
(\Q^n \oplus \{0\}) \cong \rec(P)$. 
  \end{definition}

  Given $\pi\colon \Q^n\rightarrow \Q^d$, we write
  $\widetilde{\pi}\colon \Q^n\oplus \Q\rightarrow \Q^d\oplus \Q$ for
  the map sending $({\bf x},\lambda)$ to $(\pi({\bf x}),\lambda)$.  The first
  projection $p_1\colon \Q^{n}\oplus \Q\rightarrow \Q^n$ fits into
  the diagram
 \begin{equation}
 \label{eqn:liftdiagram}
 \begin{CD}
  0 @>>> \ker(\widetilde{\pi}) @>\widetilde{\iota}>> \Q^{n}\oplus \Q
  @>\widetilde{\pi}>> \Q^{d}\oplus \Q @>>> 0 \\
  @. @V\cong VV @Vp_1VV @VVV  @. \\
  0 @>>> \ker(\pi) @>\iota>> \Q^{n} @>\pi>> \Q^d @>>> 0 \\
 \end{CD},
 \end{equation} 
 and canonically identifies $\ker(\pi)$ with $\ker(\widetilde{\pi})$.
 Write $\pi^\vee \colon (\Q^n)^*\rightarrow \ker(\pi)^*$ and
 $\widetilde{\pi}^\vee \colon (\Q^{n+1})^*\rightarrow
 \ker(\widetilde{\pi})^*$ for the maps obtained by pullback via
 $\iota$ and $\widetilde{\iota}$ respectively.  These maps fit into
 the dual diagram
 \begin{equation}
 \label{eqn:dualliftdiagram}
 \begin{CD}
   0 @<<< \ker(\widetilde{\pi})^* @<\widetilde{\pi}^\vee<< (\Q^n\oplus\Q)^* @<<< (\Q^{d}\oplus \Q)^* @<<< 0 \\
   @. @Ap_1^* A\cong A @Ap_1^* AA @AAA  @. \\
   0 @<<< \ker(\pi)^* @<\pi^\vee<< (\Q^{n})^* @<<< (\Q^d)^* @<<< 0 \\
 \end{CD}.
 \end{equation}
 Since $P\subseteq \Q^n$ and $\widetilde{P}\subseteq \Q^n\oplus\Q$,
 the normal fans satisfy $\mathcal{N}(P)\subseteq (\Q^n)^*$ and
 $\mathcal{N}\big{(}\widetilde{P}\big{)}\subseteq (\Q^n\oplus\Q)^*$
 respectively.

 \begin{proposition}
   After identifying $\ker(\pi)^*$ with $\ker(\widetilde{\pi})^*$ via
   the vertical map $p_1^*$ from diagram (\ref{eqn:dualliftdiagram}),
   the following fans coincide:
 \begin{enumerate}
 \item the fiber fan $\mathcal{N}(P,\pi(P))$;
 \item the fiber fan $\mathcal{N}(\widetilde{P},\widetilde{\pi}(\widetilde{P}))$;
 \item the fan induced by
   $\widetilde{\pi}^\vee\big{(}\mathcal{N}\big{(}\widetilde{P})\big{)}$;
 \item the fan induced by $\pi^\vee(\mathcal{N}(P))$.
 \end{enumerate}
 \end{proposition}

\begin{proof}
  We first show that fans (1) and (2) coincide.  For ${\bf q} \in
  \mathbb Q^{d}$, set $\widetilde{\bf q}:= ({\bf q},1)\in \Q^d\oplus
  \Q$.  Note that $(\widetilde{P})_{\widetilde{\bf q}} = P_{\bf
    q}$, since $P$ is the intersection of $\widetilde{P}$ with
  $\Q^d\oplus \{1\}$.  Since $\widetilde{P}$ is a cone, the polyhedron
  $(\widetilde{P})_{({\bf q},j)}$ is normally equivalent to
  $(\widetilde{P})_{({\bf q}/j,1)}$ for $j>0$.  It follows that the
  fan $\mathcal N (P,\pi(P))$ is equal to the common refinement of
  (the images under $p_1^*$ of) the normal fans of fibers
  $(\widetilde{P})_{({\bf q},j)}$ for $j>0$ and ${\bf q} \in \mathbb
  Q^d$.  Since $(\widetilde{P})_{({\bf q},j)} = \emptyset$ for $j<0$,
  the statement that the fiber fans (1) and (2) coincide follows once
  we show that $\mathcal N (P , \pi(P))$ is a refinement of the normal
  fan of $(\widetilde{P})_{({\bf q},0)}$ for all $({\bf q},0) \in
  \widetilde{\pi}(\widetilde{P})$.  Since $(\widetilde{P})_{({\bf
      q},0)}=\rec(P)_{\bf q}$, we have $({\bf q},0) \in
  \widetilde{\pi}(\widetilde{P})$ if and only if ${\bf q} \in
  \pi(\rec(P))$.  Thus, it is enough to show that $\mathcal N (P ,
  \pi(P))$ refines $\mathcal N(\rec(P)_{\bf q})$ for all ${\bf q} \in
  \pi(\rec(P))$.  We show below that for all ${\bf q} \in
  \pi(\rec(P))$ there is $\lambda > 0$ such that $\lambda{\bf q} \in
  \pi(P)$, and $\mathcal N(P_{\lambda{\bf q}})$ refines $\mathcal
  N(\rec(P)_{\lambda{\bf q}})$.  Since these normal fans equal the
  corresponding ones with $\lambda{\bf q}$ replaced by ${\bf q}$, this
  shows that $\mathcal N(P_{{\bf q}})$ and hence $\mathcal N(P,
  \pi(P))$ is a refinement of the normal fan of $\rec(P)_{\bf q}$.

  Fix ${\bf q} \in \pi(\rec(P))$. The condition that $\lambda
  \mathbf{q} \in \pi(P)$ is satisfied for all $\lambda \gg 0$, or all
  $\lambda $ if $\mathbf{0} \in P$.  Write $P=\{ \mathbf{x} \in
  \mathbb Q^n : A\mathbf{x} \geq \mathbf{b} \}$, where $A$ is an $r
  \times n$ matrix, and $\mathbf{b} \in \mathbb Q^r$.  Write $A$ in
  block form as $A=(A'|A'')$, where $A'$ is a $r \times (n-d)$ matrix,
  and $A''$ is a $r \times d$ matrix.  Then $P_{\lambda \mathbf{q}}=\{
  \mathbf{x} \in \mathbb Q^{n-d} : A'\mathbf{x} \geq \mathbf{b} -
  A''\lambda \mathbf{q}\}$, while $\rec(P)_{\lambda \mathbf{q}}=\{
  \mathbf{x} \in \mathbb Q^{n-d} : A'\mathbf{x} \geq - A''\lambda
  \mathbf{q}\}$.  We first consider the case that $A''\mathbf{q} \neq
  \mathbf{0}$.  Since $\lambda \mathbf{q} \in \pi(P)$, $P_{\lambda
    \mathbf{q}}$ and $\rec(P)_{\lambda \mathbf{q}}$ are nonempty for
  all $\lambda$, so $\mathbf{b}-A''\lambda \mathbf{q}$ and
  $-A''\lambda \mathbf{q}$ lie in some chambers of the chamber complex
  of the Gale dual of $A'$ (see \cite{BFS}).  For most generic choices
  of $\mathbf{q}$ the vector $-A''\mathbf{q}$ lies in the interior of
  a top-dimensional chamber, so for $\lambda \gg 0$ the vectors
  $-A''\mathbf{q}$ and $\mathbf{b}-A''\lambda \mathbf{q}$ lie in the
  same chamber, and thus $P_{\lambda \mathbf{q}}$ and
  $\rec(P)_{\lambda \mathbf{q}}$ are normally equivalent.  If
  $\mathbf{q}$ is not generic, then for $\lambda \gg 0$ the vector
  $\mathbf{b}-A''\lambda \mathbf{q}$ either lies in the same chamber
  as $-A''\lambda \mathbf{q}$, or in a larger dimensional chamber that
  has the chamber of $-A''\lambda \mathbf{q}$ as a face.  In the first
  case $P_{\lambda \mathbf{q}}$ and $\rec(P)_{\lambda \mathbf{q}}$ are
  normally equivalent, while in the second the normal fan of
  $P_{\lambda \mathbf{q}}$ is a refinement of that of
  $\rec(P)_{\lambda \mathbf{q}}$.  Finally, if $A''\mathbf{q} =
  \mathbf{0}$, then $\rec(P)_{\lambda \mathbf{q}}=\rec(P_{\lambda
    \mathbf{q}})$ is the recession cone of $P_{\lambda \mathbf{q}}$,
  so the refinement is automatic in this case.

  To see that fan (2) equals fan (3), apply
  Corollary~\ref{cor:coneduality} with $C=\widetilde{P}$ and
  $\widetilde{\pi}$ in place of $\pi$.  Thus, it remains to show that
  fans (3) and (4) coincide.  The proof is a simple trick.  The normal
  fan $\mathcal{N}(P)$ equals the fiber fan $\mathcal{N}(P,g(P))$ of
  the trivial linear map $g\colon \Q^n\rightarrow \Q^0$, for which the
  lift $\widetilde{g}$ is the second projection $p_2\colon \Q^n\oplus
  \Q\rightarrow \Q$.  The equality of the fans (1) and (2) established
  earlier in this proof holds for any surjective linear map, and in
  particular it holds the trivial map $g$. Thus, $p_1^*$ identifies
  $\mathcal{N}(P) = \mathcal{N}(P,g(P))$ with
  $\mathcal{N}(\widetilde{P},p_2(\widetilde{P}))$.  The equality of
  fans (2) and (3) then implies that $p_1^*$ identifies
  $\mathcal{N}(P)$ with $p_2^\vee(\mathcal{N}(\widetilde{P}))$.
  Finally, apply $\pi^\vee$ to see that $p_1^*$ identifies
  $\pi^\vee(\mathcal{N}(P))$ with $(\pi^\vee \circ
  p_2^\vee)\big{(}\mathcal{N}(\widetilde{P})\big{)}$.  If the
  upward-pointing arrow in diagram (\ref{eqn:dualliftdiagram})
  representing $p_1^*\colon (\Q^n)^*\rightarrow (\Q^n\oplus \Q)^*$ is
  replaced by the downward-pointing arrow representing $p_2^\vee\colon
  (\Q^n\oplus \Q)^*\rightarrow (\Q^n)^*$, the resulting diagram
  commutes giving $\widetilde{\pi}^\vee = p_1^*\circ \pi^\vee \circ
  p_2^\vee$.  This proves that $p_1^*$ identifies fans (3) and (4) as
  required.
 \end{proof}

 \begin{lemma}
 \label{l:heightone}
 Let $\pi\colon P\rightarrow Q$ be a linear projection of polyhedra.
 Then $\widetilde{\pi}(\widetilde{P}) = \widetilde{Q}$, so
 $\widetilde{\pi}(\widetilde{P})$ is the cone over $Q$.  Moreover, the
 fan induced by $\widetilde{\pi}(\widetilde{P})$ is the cone over the
 subdivision of $Q$ induced by $\pi(P)$, so it is obtained by taking
 the cone over each cell in that subdivision.
 \end{lemma}
 \begin{proof}
   It is easy to see that $\pi(\rec(P)) = \rec(\pi(P)) = \rec(Q)$, so
   the equality of sets $\widetilde{\pi}\big{(}\widetilde{P}\big{)} =
   \widetilde{Q}$ follows from the fact that $\widetilde{\pi}(\lambda
   {\bf x},\lambda) = (\lambda \pi({\bf x}),\lambda)$ for ${\bf x}\in
   P$. For the second statement, note first that the faces of
   $\widetilde{P}$ divide naturally into two disjoint sets: those
   arising from faces of the recession cone of $P$; and those of the
   form $\widetilde{E}$ for some face $E\subseteq P$.  The faces in
   the former set lie in $\Q^n\oplus \{0\}\subset \Q^n\oplus \Q$, so
   it is enough to show that the given decomposition of $Q$ is the
   intersection of the slice $\Q^d\oplus \{1\}\subset \Q^d\oplus \Q$
   with the fan obtained as the common refinement of the cones
   $\big{\{}\widetilde{\pi}(\widetilde{E})\subseteq \widetilde{Q} :
   E\subseteq P \text{ is a face}\big{\}}$.  The result follows since
   for every face $E\subseteq P$, the intersection of $\Q^d\oplus
   \{1\}$ with the cone $\widetilde{\pi}(\widetilde{E})\subseteq
   \widetilde{Q}$ is the cell $\pi(E)\subseteq Q$.
 \end{proof}

 \begin{proposition}
 \label{p:PveeFiberfan}
 Let $\pi\colon P\rightarrow Q$ be a linear projection of polyhedra,
 and set $C:= (\widetilde{P})^{\vee}$.  Then the fiber fan
 $\mathcal{N}( C, \widetilde{\pi}^{\vee}(C))$ is the fan induced by
 $\widetilde{\pi}(\widetilde{P})$, and thus is the cone over the
 subdivision of $Q$ induced by $\pi(P)$.
 \end{proposition}
 \begin{proof}
   Applying Corollary~\ref{cor:coneduality} to the map
   $\widetilde{\pi}^\vee\colon C\rightarrow \widetilde{\pi}^\vee(C)$
   gives $\mathcal{N}( C, \widetilde{\pi}^{\vee}(C)) =
   (\widetilde{\pi}^\vee)^\vee(C^\vee)$.  Since $\widetilde{P}$ is a
   polyhedral cone, we have $C^\vee = ((\widetilde{P})^\vee)^\vee =
   \widetilde{P}$.  The first statement follows from the fact that
   $(\widetilde{\pi}^\vee)^\vee = \widetilde{\pi}$.  The second statement
   follows from Lemma~\ref{l:heightone}.
 \end{proof}


 \section{On quotients of toric varieties by subtori}
 \label{s:toricquotients}

 \subsection{GIT quotients of quasiprojective toric varieties}
 \label{sec:toricGIT}
 We first recall the polyhedral construction of quasiprojective toric
 varieties that are projective over affine toric varieties (see
 Thaddeus~\cite[\S2.8]{Thaddeus}). We work over an algebraically closed field
 $\K$.

 Let $P\subseteq \Q^n$ be an $n$-dimensional lattice polyhedron and
 let $\widetilde{P}$ be the closure in $\Q^{n}\oplus \Q$ of the cone
 over $P$ as given in Definition~\ref{d:tildeP}.  The semigroup
 algebra $R:= \K[\widetilde{P}\cap (\Z^{n}\oplus \Z)]$ defines the
 affine toric variety $\Spec R$.  The second projection determines a
 grading of the ring $R$. Let $R_j$ be the $\K$-vector subspace
 spanned by the elements of $R$ of degree $j\in \N$. Then $X_P :=
 \Proj \bigoplus_{j\geq 0} R_j$ is the $n$-dimensional quasiprojective
 toric variety whose fan is $\mathcal{N}(P)$, the inner normal fan of
 the polyhedron $P$.  Note that $X_{\widetilde{P}}\cong \Spec R$ is
 the affine cone over $X_P$.  Also, since $R_{0}$ is isomorphic to the
 semigroup algebra $\K[\rec(P)\cap \Z^n]$ determined by the recession
 cone $\rec(P)\subseteq \Q^n$, we see that $X_P$ is projective over
 the affine toric variety $\Spec R_0$.  Thus, $X_P$ is projective if $P$ is a polytope.

 \begin{remark}
   Any quasiprojective toric variety $X$ that is projective over an
   affine toric variety arises in this way once a relatively ample
   torus-invariant divisor $D$ is chosen on $X$, in which case $R_j
   \cong H^0(X,\owe_X(jD))$ for all $j\in \N$.  
 \end{remark}

 Let $(\K^*)^n = \Spec \K[\Z^n]$ denote the dense torus of $X_P$, and
 let $T_L:= \Spec \K[L]$ be a $d$-dimensional algebraic subtorus of
 $(\K^*)^n$.  Choose generators for $L \cong \Z^d$ and write
 $\pi_{\Z}\colon \Z^n \rightarrow \Z^d$ for the map of character
 lattices induced by the inclusion $T_L\hookrightarrow (\K^*)^n$.
 Composing $\pi_{\Z}$ with the $\Z^n$-grading of the ring R arising
 from the action of $(\K^*)^n$ on $X_P$ determines a $\Z^d$-grading of
 $R$.  Given a character ${\mathbf{v}}\in \Z^d$ of the torus $T_L$,
 the GIT quotient of $X_P$ by the action of $T_L$ linearized by
 $\mathbf{v}$ is the quotient of the affine variety $\Spec(R)$ by the
 lift of $T_L$.  Specifically, begin by defining $\widetilde{\pi_{\Z}}
 \colon \Z^{n}\oplus \Z\rightarrow \Z^{d}\oplus \Z$ by sending
 $(\mathbf{x},\lambda)$ to $(\pi_{\Z}(\mathbf{x}),\lambda)$. This
 lattice map induces an inclusion of the $(d+1)$-dimensional torus
 $T_{\widetilde{L}} :=\Spec \K[(\widetilde{L}]$ into the dense torus
 $(\K^*)^{n+1} = \Spec \K[\Z^n\oplus\Z]$ of $X_{\widetilde{P}}$, where
 $\widetilde{L}= L \oplus \Z$. Characters ${\bf v}\in \Z^d$ of $T_L$
 give rise to characters $\widetilde{{\bf v}}:= ({\bf v},1)\in
 \Z^{d}\oplus\Z$ of $T_{\widetilde{L}}$.  Let $R_{j\widetilde{{\bf
 v}}}$ denote the $j\widetilde{{\bf v}}$-graded piece of the ring $R$
 with respect to the $\Z^{d}\oplus\Z$-grading of $R$.  Then
 \begin{equation}
 \label{eqn:GITdef}
 X_P \git_{\mathbf{v}} T_L=  X_{\widetilde{P}} \git_{\widetilde{\bf v}}T_{\widetilde{L}}:=
  \textstyle{\Proj \bigoplus_{j\geq 0} R_{j{\widetilde{\bf v}}}}.
 \end{equation}

 We now adopt the notation of Section~\ref{sec:GeneralDuality}.  Thus, we
 write $\pi\colon \Q^n\rightarrow \Q^d$ and $\widetilde{\pi}\colon
 \Q^{n}\oplus\Q\rightarrow \Q^{d}\oplus\Q$ for the linear maps
 obtained from $\pi_{\Z}$ and $\widetilde{\pi_\Z}$ respectively by
 extending scalars.  The following result is due to
 Kapranov--Sturmfels--Zelevinsky~\cite{KSZ} when $P$ is a polytope,
 and to Thaddeus~\cite{Thaddeus} in general.

 \begin{lemma}
 \label{lemma:XgitT}
 For ${\bf v}\in \pi(P)$, the GIT quotient $X_P\git_{\bf v}T_L$ is the toric variety with fan
 $\mathcal{N}(P_{\bf v})$.
\end{lemma}
 \begin{proof}
   We know that $X_P\git_{\bf v}T\cong
   X_{\widetilde{P}}\git_{\widetilde{\bf v}}T_{\widetilde{L}}$ from
   above.  The left-most vertical arrow in the diagram
   (\ref{eqn:liftdiagram}) sends the polyhedron
   $\widetilde{P}_{\widetilde{\bf v}}$ isomorphically onto $P_{\bf
     v}$, so the left-most vertical map in the dual diagram
   (\ref{eqn:dualliftdiagram}) identifies the normal fans
   $\mathcal{N}(\widetilde{P}_{\widetilde{\bf v}})$ and
   $\mathcal{N}(P_{\bf v})$.  The result follows from the fact that
   for $j\in \N$, the graded piece $R_{j{\widetilde{\bf v}}}$ is isomorphic to the
   $\K$-vector space spanned by the lattice points of the slice $\widetilde{P} \cap
   \widetilde{\pi}^{-1}(j\widetilde{{\bf v}})$.
\end{proof}

\begin{proposition}
  For ${\bf v}\in \pi(P)$, the fan $\mathcal{N}(P_{\bf v})$ of
  $X_{P}\git_{\bf v}T$ is isomorphic to the $\widetilde{{\bf v}}$-induced
  $\widetilde{\pi}^{\vee}$-coherent subdivision of the cone
  $\widetilde{\pi}^\vee(\widetilde{P}^{\vee})\subseteq
  \ker(\widetilde{\pi})^*$.
 \end{proposition}
 \begin{proof} 
   Corollary~\ref{cor:normalfan} shows that
   $\mathcal{N}(\widetilde{P}_{\widetilde{{\bf v}}})$ is the
   $\widetilde{{\bf v}}$-induced $\widetilde{\pi}^\vee$-coherent
   subdivision of the cone
   $\widetilde{\pi}^\vee(\widetilde{P}^{\vee})$. The result follows
   since $\mathcal{N}(P_{\bf v})$ is isomorphic to
   $\mathcal{N}(\widetilde{P}_{\widetilde{\bf v}})$.
 \end{proof}

 \begin{definition} \label{d:GITchamber} A character ${\bf v}\in
   \pi(P)$ is \emph{generic} if every ${\bf v}$-semistable point of
   $X_P$ is ${\bf v}$-stable in the sense of GIT (see
   Dolgachev--Hu~\cite[(0.2.2)]{DolgachevHu}).  For generic characters
   ${\bf v}, {\bf v}'\in \pi(P)$, we set $\mathbf{v} \sim \mathbf{v}'$
   if every ${\bf v}$-stable point of $X_P$ is ${\bf v}'$-stable and
   vice-versa.  This equivalence relation gives a polyhedral
   decomposition of $\pi(P)$, called the {\em GIT chamber
     decomposition} associated to the action of $T_L$ on $X_P$.
\end{definition}

 \begin{remark}  
 \label{r:GITchamber}
 Thaddeus~\cite{Thaddeus} showed that the GIT chamber decomposition is
 the subdivision of $Q$ induced by $\pi(P)$ in the sense of
 Definition~\ref{d:subdivisioninduced}. For analogous results
 for general GIT quotients, see Thaddeus~\cite{ThaddeusGeneral} or
 Dolgachev--Hu~\cite{DolgachevHu}.
\end{remark}

 \subsection{The toric Chow quotient}
 \label{sec:toricChow}
 In this section we construct the toric Chow quotient of a
 quasiprojective toric variety by a subtorus.  We first introduce a
 subcategory of the category of $\K$-schemes.  Fix $T := (\K^*)^n$ for
 $n\in \N$.

 \begin{definition}
   Consider the category $\mathcal{C}_T$ whose objects are Noetherian
   $\K$-schemes $X$ with a $T$-action $\sigma_X \colon T\times
   X\rightarrow X$ and an irreducible $T$-invariant component
   $X_0\subseteq X$, such that the restriction of the $T$-action to
   $X_0$, denoted $\sigma_{X_0} \colon T\times X_0\rightarrow X_0$,
   gives $X_0$ the structure of a not-necessarily normal toric variety
   with dense torus $T$.  The morphisms of $\mathcal{C}_T$ are
   $T$-equivariant proper morphisms $f\colon X\rightarrow X'$ over
   $\K$ that induce a birational morphism $f\vert_{X_0}\colon
   X_0\rightarrow X'_0$ between the toric components.
 \end{definition}

 \begin{lemma} \label{l:fiberprodsexist}
 Fiber products exist in the category $\mathcal{C}_T$. Moreover, the
 fiber product in $\mathcal{C}_T$ coincides with the fiber product in
 the category of $\K$-schemes.
 \end{lemma}
 \begin{proof}
 For an object $S$ in $\mathcal{C}_T$, let $f\colon X\rightarrow S$
 and $f'\colon X'\rightarrow S$ be morphisms in $\mathcal{C}_T$.  We
 must show that the fiber product in the category of $\K$-schemes
 $Z:=X\times_S X'$ is an object of $\mathcal{C}_T$, and that the
 canonical projections $p\colon Z\rightarrow X$ and $q\colon
 Z\rightarrow X'$ are morphisms in $\mathcal{C}_T$.

 To construct the component $Z_0\subseteq Z$, note that the schemes
 $X$, $X'$ and $S$ contain open subvarieties $T_X$, $T_{X'}$ and $T_S$
 respectively, each isomorphic to the torus $T$. The fiber product
 $T_Z:= T_X\times_{T_S} T_{X'}$ is canonically isomorphic to the
 subscheme $p^{-1}(T_X)\cap q^{-1}(T_{X'})$ of $Z$ by \cite[Chap.\ I,
 Coro 3.2.3]{EGA}, so it is open in $Z$. Furthermore, the
 $T$-equivariant birational morphisms $f\vert_{X_0}\colon
 X_0\rightarrow S_0$ and $f'\vert_{X'_0}\colon X'_0\rightarrow S_0$
 restrict to give isomorphisms $T_X\cong T_S$ and $T_{X'}\cong T_S$,
 so $T_Z$ is isomorphic to $T\times_T T\cong T$ and hence is
 irreducible.  In particular, the $n$-torus $T_Z\subseteq Z$ is dense
 in some component of $Z$ that we denote $Z_0$.

 We next show that $Z_0$ is reduced.  Since this is a local question,
 and $T_Z \subseteq Z_0$ is reduced, we can reduce to the case where
 $Z_0=\Spec(A)$ is an irreducible Noetherian $\K$-scheme with a
 nonempty open subscheme $W$ for which the induced scheme structure on
 $W$ is reduced.  We may assume that $W=\Spec(A_f)$ for some $f \in
 A$.  Writing $A=\K[x_1,\dots,x_m]/I$ for some $I$, we note that $I$
 is primary since $Z_0$ is irreducible, and since $W$ is nonempty $f
 \not \in \sqrt I$.  Suppose $Z_0$ is not reduced, so there exists $g
 \not \in I$ with $g^k \in I$ for some $I$.  Then $(g/1)^k \in I_f$,
 so since $W$ is reduced we have $g/1 \in I_f$.  But then $g/1=y/f^l$
 for some $l$ and some $y\in I$.  This means that $gf^l \in I$, which
 contradicts $I$ being primary, $g \not \in I$, and $f \not \in \sqrt
 I$.  We thus conclude that $Z_0$ is reduced.

 We now describe the torus action.  The fiber product $Z$ embeds as a
 closed subscheme of the product $X\times_{\K} X'$.  Since the
 morphisms $f$ and $f'$ are $T$-equivariant, it can be shown that the
 product $T$-action on $X\times_{\K} X'$ restricts to give an action
 $\sigma_Z\colon T\times Z\rightarrow Z$.  This restricts to give an
 action $T\times T_Z\rightarrow T_Z$ that coincides with the
 multiplicative structure on the torus $T_Z$ after identifying
 $T_Z\cong T$. Thus, $\sigma_Z\colon T\times Z\rightarrow Z$ extends
 the natural multiplicative structure of $T_Z$.  To see that $Z_0$ is
 $T$-invariant we argue by continuity as follows.  The component
 $Z_0\subseteq Z$ is closed and contains $T_Z$, so $\sigma_Z^{-1}(Z)$
 is a closed subscheme of $T\times Z$ that contains
 $\sigma_Z^{-1}(T_Z) = T\times T_Z$. In particular, $\sigma_Z^{-1}(Z)$
 contains the closure $T\times Z_0$ of $T\times T_Z$, so the image of
 $T\times Z_0$ under $\sigma_Z$ lies in $Z_0$.  The resulting action
 $\sigma_{Z_0}\colon T\times Z_0\rightarrow Z_0$ extends the
 multiplicative structure of $T_Z$ on itself by the above, so $Z_0$ is
 a not-necessarily-normal toric variety.

 It remain to prove that the projections $p$ and $q$ are morphisms in
 $\mathcal{C}_T$.  Properness of $p$ and $q$ follows from properness
 of $f$ and $f'$ by base extension.  The $T$-action on $Z$ was
 constructed to ensure that $p$ and $q$ are $T$-equivariant.  Their
 restrictions $p\vert_{T_Z} \colon T_Z \rightarrow T_X$ and
 $p'\vert_{T_Z} \colon T_Z \rightarrow T_{X'}$ are isomorphisms, so
 both $p$ and $q$ are birational on the toric components as required.
 \end{proof}

 We now return to the GIT set-up from Section~\ref{sec:toricGIT},
 where $\pi\colon P\rightarrow Q$ is a linear projection of polyhedra,
 with $L \cong \mathbb Z^d$ the image of the corresponding map
 $\pi_{\mathbb Z} \colon \mathbb Z^n \rightarrow \mathbb Z^d$.  Set
 $M:=\ker(\pi_{\Z})$.  For ${\bf v}\in Q$, the GIT quotient
 $X_P\git_{\bf v}T_L = X_{P_{\bf v}}$ is the toric variety with dense
 torus $T_M$ defined by the fan $\mathcal{N}(P_{\bf v})$.
 
 Let $\poset$ be the face poset of the subdivision of $Q$ induced by
 $\pi(P)$, with the faces on the boundary of $Q$ removed.  Thus $\tau
 \prec \sigma$ if $\tau$ is a face of $\sigma$. To each $\sigma \in \poset$
 we associate the toric variety $X_{\sigma}:=X_{P_{\bf v}}$ for any
 $\mathbf{v} \in \sigma$.  If $\tau$ is a face of $\sigma$, then there
 is a proper, birational toric morphism from $X_{\sigma}$ to
 $X_{\tau}$ by Thaddeus~\cite[Theorem 3.11 and
 Corollary~3.12]{Thaddeus}.  Therefore $\poset$ gives a directed
 system of $\K$-schemes.  We are interested in the inverse limit of
 this system in the category of $\K$-schemes.

 \begin{proposition}
   \label{prop:invlimit} 
   Let $T_M$ be the dense torus in $X_{\sigma}$ for $\sigma \in
   \poset$. Then the inverse limit in the category of $\K$-schemes
 \begin{equation}
 \label{eqn:invlimit}
 Z:= \varprojlim_{\sigma\in \poset}
   X_{\sigma}
 \end{equation}
 exists and is an object of the category $\mathcal{C}_{T_M}$. In
 particular, $Z$ has an irreducible component $Z_0$
 that is a not-necessarily-normal toric variety with dense torus
 $T_M$.
 \end{proposition}
 \begin{proof}
   For each $\sigma \in \poset$, the toric variety $X_{\sigma}$ is an
   object of the category $\mathcal{C}_{T_M}$.  The proof is by
   induction on the number of elements of $\poset$.  The base case is
   when $\poset$ has only one maximal element $\sigma$.  Then
   $Z=X_{\sigma}$ exists and is an object of $\mathcal{C}_{T_M}$.
   Suppose now that the inverse limit exists, and is an object in
   $\mathcal{C}_{T_M}$, for any such poset with fewer elements than
   $\poset$ whose elements are toric varieties with dense torus $T_M$,
   and whose maps are proper birational toric morphisms.  Fix a
   maximal element $\sigma$ of $\poset$, and let $\poset'=\poset
   \setminus \sigma$, and let $\poset''$ be the subposet of all $\tau
   \in \poset$ with $\tau \prec \sigma$.  Let $Z_{\poset'}$ and
   $Z_{\poset''}$ be the corresponding inverse limits.  Define 
   $$Z':=Z_{\poset'} \times_{Z_{\poset''}} X_{\sigma}.$$
   
   Lemma~\ref{l:fiberprodsexist} and induction show that $Z'$ is an
   object in $\mathcal C_{T_M}$.  It remains to show that $Z'$ is the
   inverse limit~(\ref{eqn:invlimit}) in the category of $\K$-schemes.
   Let $Y$ be any $\K$-scheme with maps $f_{\tau}$ to each $X_{\tau}$
   for $\tau \in \poset$ that commute appropriately with the maps in
   $\poset$.  By the universal property of the inverse limits
   $Z_{\poset'}$ and $Z_{\poset''}$ we get a unique map from
   $Z_{\poset'}$ to $Z_{\poset''}$ that commutes with the maps in
   $\poset$, and a unique map from $Y$ to $Z_{\poset'}$ and
   $Z_{\poset''}$ whose compositions with the maps from $Z_{\poset'}$
   to $Z_{\poset''}$ and from $X_{\sigma}$ to $Z_{\poset''}$ agree.
   Then the universal property of fiber products gives a unique map
   $\phi$ from $Y$ to $Z'$ such that the composition of $\phi$ with the
   maps from $Z'$ to each $X_{\tau}$ equals the corresponding
   $f_{\tau}$, so $Z'$ is the inverse limit.
\end{proof}

 \begin{definition}
 \label{defn:toricchow}
 The \emph{toric Chow quotient} of the toric variety $X_P$ by the
 subtorus $T_L$ is the not-necessarily-normal toric variety $X_P/T_L:=
 Z_0$ arising as the toric component of the inverse limit $Z$ as in
Proposition~\ref{prop:invlimit}.
 \end{definition}

 Corollary~\ref{coro:chow} below justifies the choice of terminology; see also \cite{HaimanSturmfels}.
 First, we explain the link between the toric Chow quotient and fiber
 fans.

 \begin{proposition}
 \label{p:toricChow}
 Let $\overline{X_P/T_L}$ be the normalization of the toric Chow
 quotient.  Then $\overline{X_P/T_L}$ is the toric variety with fan
 $\mathcal{N}(P,\pi(P))$.
 \end{proposition}
 \begin{proof}
   Write $\overline{Z_0} = \overline{X_P/T_L}$.  For each $\sigma \in
   \poset$ there is a proper, birational toric morphism
   $\overline{Z_0}\rightarrow X_{\sigma}$, so the fan of
   $\overline{Z_0}$ refines the fan of $X_{\sigma}$ by
   \cite[Proposition 2.4]{Fulton}.  In particular, the fan of
   $\overline{Z_0}$ refines the common refinement
   $\mathcal{N}(P,\pi(P))$ of these fans.
   
   The proposition follows once we show that $\mathcal{N}(P,\pi(P))$
   refines the fan of $\overline{Z_0}$.  It is enough to construct a
   proper, birational toric morphism $X_{CR}\rightarrow \overline{Z_0}$,
   where $X_{CR}$ is the toric variety with fan $\mathcal{N}(P,\pi(P))$.
   Since $\mathcal{N}(P,\pi(P))$ refines the fan of each $X_{\sigma}$, there is a proper, birational toric morphism
   $X_{CR}\rightarrow X_{\sigma}$.  The universal property of the
   inverse limit $Z$ in the category $\mathcal{C}_{T_M}$ then gives a
   proper $T_M$-equivariant morphism $\phi\colon X_{CR}\rightarrow Z$
   that is birational onto $Z_0$.  Since $X_{CR}$ is normal, the
   universal property of normalization gives a morphism
   $\overline{\phi}\colon X_{CR}\rightarrow \overline{Z_0}$ that is both
   proper (since $\phi$ is proper and the normalization is separated)
   and birational (since $\phi$ is birational onto $Z_0$ and the
   normalization is finite).  It remains to show that
   $\overline{\phi}$ is a toric morphism.

   Since $\phi$ is a dominant birational morphism of
   not-necessarily-normal toric varieties, it is given locally by an
   inclusion of subsemigroups $S\rightarrow S'$ of $M$, where $\Spec
   \K[S]$ is a chart on $Z_0$, $\Spec \K[S']$ is a chart on $X_{CR}$ and $T_M
   = \Spec \K[M]$.  Note that $S'$ is saturated since $X_{CR}$ is normal.
   The normalization $\overline{Z_0}\rightarrow Z_0$ is given locally by
   the inclusion of $S$ in its saturation $S_{\sat}$. Since $S_{\sat}$
   is the smallest saturated subsemigroup of $M$, the inclusion
   $S\hookrightarrow S'$ factors through $S_{\sat}$. The induced
   semigroup morphism $S_{\sat}\rightarrow S'$ gives the local
   description of the induced map $\overline{\phi}\colon X_{CR} \rightarrow
   \overline{Z_0}$, so $\overline{\phi}$ is a toric morphism.
 \end{proof}

 \begin{corollary}
 \label{coro:chow}
 Let $P$ be a polytope. Then the toric Chow quotient $X_P / T_L$ is equal to
 the Chow quotient in the sense of Kapranov~\cite[Section 0.1]{Kapranov}.
 \end{corollary}
 \begin{proof}
   Write $\Chow$ for the Chow quotient.  For a polytope $P$,
   Kapranov-Sturmfels-Zelevinsky~\cite[Proposition~2.3]{KSZ} proved
   that the normalization of the Chow quotient is the toric variety
   defined by the fan $\mathcal{N}(P,\pi(P))$, so $X_P / T_L$ and
   $\Chow$ share the same normalization. To see that they coincide,
   recall from \cite[Corollary~4.3]{KSZ} that $\Chow$ is a subvariety
   of the inverse limit $Z$ defined in equation (\ref{eqn:invlimit}).
   Since $\Chow$ is irreducible and contains the dense torus $T_Z$, it
   is a subvariety of the component $Z_0 = X_P/T_L$ of $Z$. If this
   inclusion $\Chow\subseteq X_P/T_L$ were strict, $\Chow$ and
   $X_P/T_L$ could not share the same normalization, hence $\Chow =
   X_P / T_L$ as claimed.
 \end{proof}

 \begin{remark}
   The statement of Corollary~\ref{coro:chow} was known to the authors
   of \cite[Section~4]{KSZ}, though they did not supply a proof.
Thus, for a polyhedron $P$, the toric Chow quotient $X_P/T_L$ is
the appropriate generalization of the Chow quotient $\Chow$. 
 \end{remark}

\subsection{Proof of the main theorems}
 \label{s:mainproof}

 We now prove the main theorems.  Recall that the fan of a
 not-necessarily-normal toric variety, such as the toric Chow
 quotient, is the fan of its normalization.

\begin{theorem}
 \label{thm:main}
 Let $\pi\colon
 P\rightarrow Q$ be a linear projection of polyhedra.  Set $P':=
 (\widetilde{P})^{\vee}$ and $L':= \image(\widetilde{\pi}^\vee)$.
 The cone over the GIT chamber decomposition arising from the
 $T_L$-action on $X_P$ is the fan of the toric Chow quotient
 $X_{P'}/T_{L'}$.
\end{theorem}
\begin{proof}
  The fan defining the toric Chow quotient $X_{P'}/T_{L'}$ is
  $\mathcal{N}(P',\widetilde{\pi}^\vee(P'))$ by
  Proposition~\ref{p:toricChow}.  Now Proposition~\ref{p:PveeFiberfan}
  implies that the fan defining the toric Chow quotient
  $X_{P'}/T_{L'}$ is equal to the cone over the subdivision of $Q$
  induced by $\pi(P)$.  The result follows from
  Remark~\ref{r:GITchamber}, since the subdivision of $Q$ induced by
  $\pi(P)$ is the GIT chamber decomposition arising from the
  $T_L$-action on $X_P$.
\end{proof}

The case when $P$ is a polyhedral cone (so $X_P$ is affine) is of
particular interest.  For example, Cox~\cite{Cox} showed that every
simplicial, projective toric variety can be constructed as a GIT
quotient of affine space by a subtorus.  In the affine case, Theorem~\ref{thm:main} can be strengthened as follows:
  
 \begin{theorem}[Duality for affine toric quotients]
 \label{thm:affineduality}
   Let $P\subseteq \Q^n$ be a polyhedral cone and $\pi\colon
   P\rightarrow Q$ a linear projection.  Set $P' := P^\vee$ and $L':=
   \image(\pi^\vee)$.  Then
 \begin{enumerate}
 \item[\one] the fan of the toric Chow quotient $X_{P}/T_{L}$ is equal
   to the GIT chamber decomposition arising from
   the action of $T_{L'}$ on $X_{P'}$; and
 \item[\two] the GIT chamber decomposition arising from the
   $T_L$-action on $X_P$ is the fan of the toric Chow quotient
   $X_{P'}/T_{L'}$.
 \end{enumerate}
 \end{theorem}
 \begin{proof}
   For part \one, note that the fan of
   $X_{P}/T_{L}$ is the fiber fan $\mathcal{N}(P,\pi(P))$ by
   Proposition~\ref{p:toricChow}.  Corollary~\ref{cor:coneduality}
   shows that $\mathcal{N}(P,\pi(P))$ equals the fan induced by
   $\pi^\vee(P^\vee)$.  This is precisely the GIT chamber
   decomposition arising from the action of $T_{L'}$ on
   $X_{P^\vee}$ according by Remark~\ref{r:GITchamber}.

   For part \two, note that the fan of the toric Chow quotient
   $X_{P'}/T_{L'}$ is the fiber fan
   $\mathcal{N}(P^\vee,\pi^\vee(P^\vee))$.  This equals
   $(\pi^\vee)^\vee((P^\vee)^\vee) = \pi(P)$ by applying
   Corollary~\ref{cor:coneduality} with $C = P^\vee$ and $\pi^\vee$
   playing the role of $\pi$.
 \end{proof}

\begin{remark}
  If $X$ is a nonnormal quasiprojective toric variety then one can
  give similar combinatorial descriptions of both GIT quotients by
  subtori $X \git_{\mathbf{v}}T$, and of the normalization of the
  toric Chow quotient; normality played no essential role in either
  construction.  However the description of the GIT chamber
  complex in Definition~\ref{d:GITchamber} is now a priori only a
  coarsening of the true GIT chamber complex, as the quotient toric
  varieties are no longer determined solely by their fans.
\end{remark}

\def\cprime{$'$}
\providecommand{\bysame}{\leavevmode\hbox to3em{\hrulefill}\thinspace}
\providecommand{\MR}{\relax\ifhmode\unskip\space\fi MR }
\providecommand{\MRhref}[2]{%
  \href{http://www.ams.org/mathscinet-getitem?mr=#1}{#2}
}
\providecommand{\href}[2]{#2}


\end{document}